\documentclass{amsart}%

\usepackage{graphicx}
\usepackage{latexsym}
\usepackage{amsmath}
\usepackage{amsfonts}
\usepackage{amssymb}
\usepackage{amsthm}
\usepackage{bm} 
\usepackage{url}

\input epsf.sty

\newtheorem{theorem}{Theorem}[section]
\newtheorem{proposition}[theorem]{Proposition}

\title[Optimal Packings on Flat Tori]{Optimal Packings of Two to Four Equal Circles on Any Flat Torus}

\author[Brandt]{Madeline Brandt}
\address{Department of Mathematics \\
University of California, Berkeley\\
Berkeley, CA 94720-3840}
\email[Madeline Brandt]{brandtm@berkeley.edu}

\author[Dickinson]{William Dickinson}
\address{Department of Mathematics\\
Grand Valley State University\\
Allendale, MI 49401}
\email[William Dickinson]{dickinsw@gvsu.edu}
\thanks{The authors were partially supported by National Science Foundation grants DMS-0451254 and DMS-1262342. The second author would like to thank David Austin for reading drafts of this article and writing a Python interface for using the program PiScript (written by Bill Casselman) that we used to produce the illustrations in this article. We all express our gratitude to Bob Connelly for inspiring this work. We also thank the referees for their comments that improved this work.}

\author[Ellsworth]{AnnaVictoria Ellsworth}
\address{254 Windsor St Unit 3L\\
Cambridge, MA 02139}
\email[AnnaVictoria Ellsworth]{AnnaVictoria.ellsworth@gmail.com}

\author[Kenkel]{Jennifer Kenkel}
\address{Department of Mathematics\\
University of Utah \\
Salt Lake City, Utah 84112-0090}
\email[Jennifer Kenkel]{kenkel@math.utah.edu}

\author[Smith]{Hanson Smith}
\address{Department of Mathematics\\
University of Colorado \\
Boulder, Colorado 80309-0395}
\email[Hanson Smith]{hanson.smith@colorado.edu}

\keywords{Equal circle packing, flat torus, packing graph, rigidity theory}
\subjclass[2010]{52C15} 

\date{\today}
\begin{document}

\begin{abstract}
We find explicit formulas for the radii and locations of the circles in all the optimally dense packings of two, three or four equal circles on any flat torus, defined to be the quotient of the Euclidean plane by the lattice generated by two independent vectors. We prove the optimality of the arrangements using techniques from rigidity theory and topological graph theory. 
\end{abstract}

\maketitle

\section{Introduction and Main Result}

We consider a problem from discrete geometry in the area of equal circle packing in a domain. The domain that we consider is the quotient of the Euclidean plane by a full-rank lattice, called a flat torus.  We determine all optimally dense arrangements of two, three, and four equal circles on any flat torus. Theorem~\ref{thm:MainResult} states explicit formulas for the optimal radius on any torus.  The complexity of the expressions for the optimal radii is a consequence of the explicit formulas for the coordinates of the circle centers achieving the optimal densities given in Section~\ref{sec:Existence}. 

\begin{figure}
\includegraphics{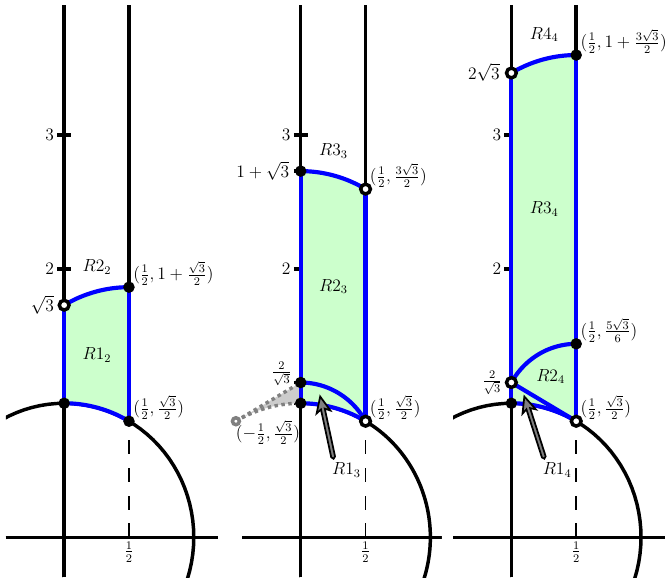}
\caption{The different expressions for the optimal radii from Theorem~\ref{thm:MainResult} break the moduli space for un-oriented tori into various regions for 2 circles (left), 3 circles (middle) and 4 circles (right). The unfilled dots on the perimeter of each moduli space indicate tori on which the optimal packing is one in which each circle is tangent to six circles (i.e. the packing lifts to the triangular close packing in the plane). 
} \label{fig:ModuliSpaces}
\end{figure}

The geometry of a flat torus influences the nature of the optimal packings it admits.  Therefore, we must first carefully describe how the optimal arrangements break the moduli space of flat tori apart.  
The moduli space of un-oriented flat tori can be represented as a strip in the $x$-$y$ plane where $x^2+y^2\geq 1$, $y>0$, and $0 \leq x \leq \frac{1}{2}$.  Within this moduli space the optimal arrangements determine regions, bounded by line and circle segments, where the optimal radii are governed by different formulas.  These regions are defined in Table~\ref{tab:regions} and pictured in Figure~\ref{fig:ModuliSpaces}. Using these regions we can state the main result of this article.

\begin{table}
\scalebox{.75}{
\renewcommand{\arraystretch}{2}
\begin{tabular}{|c|c|} 
\multicolumn{2}{c}{Two Circles} \\ \hline
Region $R1_2$ & $\sqrt{1-x^2} \leq y <  \sqrt{1-(x-\frac{1}{2})^2}+\frac{\sqrt{3}}{2}$ \\ \hline 
Region $R2_2$ & $\sqrt{1-(x-\frac{1}{2})^2}+\frac{\sqrt{3}}{2} \leq y$ \\ \hline
\multicolumn{2}{c}{Three Circles} \\ \hline
Region $R1_3$ & $ \sqrt{1-x^2} \leq y <  \sqrt{\frac{1}{3}-x^2}+\frac{\sqrt{3}}{3}$ \\ \hline
Region $R2_3$ & $\sqrt{\frac{1}{3}-x^2} +\frac{\sqrt{3}}{3} \leq y < \sqrt{1-x^2} +\sqrt{3}$\\ \hline
Region $R3_3$ & $ \sqrt{1-x^2} +\sqrt{3} \leq y$\\ \hline
\multicolumn{2}{c}{Four Circles} \\ \hline
Region $R1_4$ & $\sqrt{1-x^2} \leq y < \frac{2-x}{\sqrt{3}} $ \\ \hline
Region $R2_4$ & $\frac{2-x}{\sqrt{3}} \leq y < \sqrt{\frac{1}{3}-\left(x-\frac{1}{2}\right)^2}+\frac{\sqrt{3}}{2}$\\ \hline
Region $R3_4$ & $\sqrt{\frac{1}{3}-\left(x-\frac{1}{2}\right)^2}+\frac{\sqrt{3}}{2} \leq y < \sqrt{1-\left(x-\frac{1}{2}\right)^2}+\frac{3 \sqrt{3}}{2}$\\ \hline
Region $R4_4$ & $\sqrt{1-\left(x-\frac{1}{2}\right)^2}+\frac{3 \sqrt{3}}{2} \leq y$\\ \hline
\end{tabular}}
\vspace{.1in}
\caption{These are the equations describing the boundaries of each region in Figure~\ref{fig:ModuliSpaces}. These regions break the moduli space for un-oriented flat tori ($x^2+y^2\geq 1$, $y>0$, and $0 \leq x \leq \frac{1}{2}$) into various regions bounded by line and circle segments.  The subscript in the region name refers to the number of circles in the packing.}\label{tab:regions}
\end{table}

\newpage
 
\begin{theorem}\label{thm:MainResult} Let $\mathcal{P}_k(x,y)$ be a packing of $k$ equal circles on a flat torus that is the quotient of the plane by the vectors $\textbf{\emph{v}}_1 = \langle 1, 0\rangle$ and $ \textbf{\emph{v}}_2 =\langle x, y\rangle$ where $x^2+y^2\geq 1$, $y>0$, and $0 \leq x \leq \frac{1}{2}$. If $r_k(x,y)$ is the least upper bound of the radius then we have the following expressions for the optimal radii in different regions given explicitly in Table~\ref{tab:regions} and shown in Figure~\ref{fig:ModuliSpaces}. 
\begin{align*}
	r_2(x,y) =& 
		\begin{cases}
			\frac{\sqrt{x^2+y^2}}{4y} \sqrt{(x-1)^2+y^2}
				& \text{~\hspace{.615in}in Region $R1_2$;}\\
		         \frac{1}{2} 
				& \text{~\hspace{.615in}in Region $R2_2$.}
		\end{cases} \\
	r_3(x,y) =& 
		\begin{cases}
			\frac{\sqrt{x^2+y^2}}{2 (y+\sqrt{3} x)} \sqrt{\left(x+\frac{1}{2}\right)^2+\left(y-\frac{\sqrt{3}}{2}\right)^2}
				 &  \text{in Region $R1_3$;}\\
			\frac{1}{6} \sqrt{9 x^2+\left(y-\sqrt{3+4y^2-12 x^2}\right)^2}
				& \text{in Region $R2_3$;}\\
			\frac{1}{2} 
				& \text{in Region $R3_3$.}
		\end{cases} \\
	r_4(x,y) =& 
		\begin{cases}
			\frac{\sqrt{x^2+y^2}}{2( y- \sqrt{3} x)} \sqrt{\left(x-\frac{1}{2}\right)^2+\left(y-\frac{\sqrt{3}}{2}\right)^2}
				 &  \text{~\hspace{-0.015in}in Region $R1_4$;}\\
				 				& \\
			\frac{\sqrt{A_2-\sqrt{A_2^2-16 B_2}}}{4 \sqrt{2}} \\
				& \vspace{-0.6in} \\
				& \text{~\hspace{-0.015in}where $A_2=2 \left(y\sqrt{3}-x\right)^2-\left((x-1)\sqrt{3}+y\right)^2+3$}\\
				& \text{~\hspace{-0.015in}and $B_2=(x^2+y^2) (\left(x-\frac{3}{2}\right)^2+\left(y-\frac{\sqrt{3}}{2}\right)^2)$}\\
				& \text{~\hspace{-0.015in}in Region $R2_4$;}\\
				& \\
			\frac{\sqrt{A_3-\sqrt{A_3^2-16 B_3}}}{8 \sqrt{2}} 
			        & \vspace{-0.4in} \\
				& \text{~\hspace{-0.015in}where $A_3=9+5 y^2- (2x-1)^2$}\\
				& \text{~\hspace{-0.015in}and $B_3=\left((x-2)^2+y^2\right) \left((x+1)^2+y^2\right)$}\\
				& \text{~\hspace{-0.015in}in Region $R3_4$;}\\
				& \\
			\frac{1}{2} 
				& \text{~\hspace{-0.015in}in Region $R4_4$.}
		\end{cases}
\end{align*}
\end{theorem}

The expressions from Theorem~\ref{thm:MainResult} for the optimal radii agree with the work of Heppes~(\cite{heppes}) who presents the optimal radii\footnote{There is a typo in Heppes' formulas for the optimal radius for the rectangular tori on the edges of region $R1_3$ and $R1_4$ (in the notation of this article) for three and four circles. Each of his formulae for these radii need to be multiplied by $\frac{1}{4}$.} and arrangements for two, three and four equal circle packings on any rectangular torus (the quotient of the plane by perpendicular lattice vectors).  Heppes used the same techniques as Melissen~(\cite{melissen}) who determined the optimally dense arrangements of 1 to 4 equal circles in a square flat torus (the quotient of the plane by unit and perpendicular lattice vectors).   Both Heppes and Melissen utilized a technique similar to the one commonly used in proving the optimality of an arrangement of equal circles packed into a unit square.  This technique centers on knowing a candidate for the densest arrangement of circles in the square or rectangle which establishes a lower bound on the diameter of circles in the densest arrangement.  One next cleverly uses this diameter to partition the square or rectangle into regions with an appropriate diameter to prove the global optimality of the candidate arrangement. The approach employed in this article is fundamentally different. Here we prove the optimality of the arrangements of 2, 3 and 4 equal circles on a flat torus using techniques from rigidity theory and topological graph theory.  

Besides the work of Melissen and Heppes, packings on flat tori have been studied by other authors.  Przeworski~(\cite[Theorem 2.3]{przeworski}) determines the optimal arrangements of two equal circles on any flat torus. 
  Dickinson \textit{et al.}~(\cite{dickinson2,dickinson3}) determine the optimal packings of 1--5 equal circles on the square flat torus and 1--6 equal circles on a triangular flat torus (the quotient of the plane where the lattice vectors are unit and form an angle of $\frac{\pi}{3}$). The arrangements and radii from these articles agree with the results presented here. In addition, Musin and Nikitenko~(\cite{musin}) use similar techniques coupled with a computer algorithm to numerically determine the optimal arrangements of 6, 7 and 8 equal circles on a flat square torus. Lubachevsky \textit{et al.}~(\cite{graham}) explore packings with large numbers (50-10,000) of equal circles packed on a square torus (among other domains). They used a billiards algorithm to discover their arrangements and they discuss large scale patterns as there is little hope of proving optimality. For similar explorations from a physics point of view, see the article by Donev \textit{et al.}~(\cite{donev}). Articles~\cite{sloan,gensane} explore optimal packings of squares on the square flat torus.

The present work use similar techniques as in \cite{dickinson2}, \cite{dickinson3}, and \cite{musin} to discover and prove the optimality of equal circle packings on any flat torus. In Section~\ref{sec:def}, we review basic terms and definitions  and outline the structure of the proof of Theorem~\ref{thm:MainResult}. Section~\ref{sec:rigid} recalls results from rigidity theory and delineates the combinatorial properties of packing graphs associated to optimally dense packings whose packing graphs do not contain a loop. Section~\ref{sec:selftangent} characterizes all equal circle packings that contain at least one self-tangent circle or whose packing graphs contain at least one loop. We use topological graph theory in Section~\ref{sec:TopologicalGraphTheory} to enumerate of all 2-cell embeddings of combinatorial graphs onto a topological torus and conclude that section with tools that eliminate many of the embedded graphs from being associated with any equal circle packing on any torus.  Section~\ref{sec:nonoptimal} contains a discussion on the remaining embedded graphs that are not associated with an optimal equal circle packing (including those that are locally but not globally maximally dense).  Finally, in Section~\ref{sec:Existence} we demonstrate realizations of the remaining embedded graphs as optimal equal circle packings and thus prove the existence of packings achieving the optimal radii given in Theorem~\ref{thm:MainResult} and give the coordinates of the circle centers. 

\section{Basic Notions And An Overview Of The Proof Of Theorem~\ref{thm:MainResult}}\label{sec:def}

In this section we review some terminology and recall some basic facts about circle packings. We then outline the proof of Theorem~\ref{thm:MainResult}. 

\subsection{Basic Notions}
The quotient of the Euclidean plane by a lattice generated by two independent vectors $\textbf{v}_1$ and $\textbf{v}_2$ is called a \textbf{flat torus}.  When the vectors are perpendicular the quotient is called a \textbf{rectangular flat torus}; when they are perpendicular and unit, the quotient is called a \textbf{square flat torus}; when they are  unit and form an angle of $\frac{\pi}{3}$ the quotient is called a \textbf{triangular flat torus}.  A \textbf{fundamental domain} of a flat torus is the set of points in the Euclidean plane, $\left\{ t_1\textbf{v}_1 + t_2 \textbf{v}_2 \mid t_1, t_2 \in \mathbb{R}, 0 \leq t_1, t_2 < 1 \right\}$. To specify the location of a circle on torus, we will actually give a location in the universal cover (the Euclidean plane) which will serve as a representative of the equivalence class of the lifts of the location from the torus. In a slight abuse of nomenclature, we shall also say that two points in the Euclidean plane that are in the same equivalence class are also lifts of each other (this implies that they differ by a vector in the lattice).  Notice that we can lift any packing from a flat torus to a periodic packing in the Euclidean plane by lifting all the circles in the packing in all possible ways. This means that the density of a packing on a flat torus must be less than the maximum density of a packing of the Euclidean plane.  The L. Fejes T\'oth-Thue Theorem (\cite{toth, thue}) states that the density of all packings of equal circles in the Euclidean plane is less than the density of the triangular close packing, where each circle is tangent to six others. The triangular close packing has density $\frac{\pi}{\sqrt{12}}$ and therefore the packing density on the torus cannot exceed this bound.

The \textbf{standard basis} for a lattice is one where $\mathbf{v}_1= \langle 1,0 \rangle$ and $\mathbf{v}_2 =\langle x,y \rangle$ where $x^2+y^2 \geq 1$, $y>0$ and $-\frac{1}{2} < x \leq \frac{1}{2}$. Every oriented lattice can be transformed by scaling and the action of $SL(2,\mathbb{Z})$ to have a standard basis (see, for example, \cite[Theorem 2.3]{przeworski} or \cite[Theorem 2.7.1]{Jost}). As we are working with unoriented tori, we restrict to the region where $0 \leq x \leq \frac{1}{2}$. The optimal arrangements naturally break this restricted moduli space of flat tori into regions with different expressions for the optimal radius. See Figure~\ref{fig:ModuliSpaces}. 

For a given flat torus, an arrangement of equal circles forms a \textbf{packing} on the torus if the interiors of the circles are disjoint. The \textbf{density of a packing} is the ratio of the area of the circles to the area of the flat torus. Notice that the transformations used to put a lattice (with a periodic circle packing) into a standard form preserve the density of that packing. We define two packings with the same number of circles to be \textbf{$\bm{\epsilon}$-close} if there is a one-to-one correspondence between the circles, so that corresponding circles have centers that are all within a distance of $\epsilon$ of each other. We define a packing $\mathcal{P}$ to be optimal or \textbf{locally maximally dense} if there exists an $\epsilon > 0$ so that all $\epsilon$-close packings of equal circles have a packing density no larger than that of $\mathcal{P}$. A packing $\mathcal{Q}$ is \textbf{globally maximally dense} if it is the densest possible packing. Rather than searching directly for the globally maximally dense packings, our techniques allow us to determine all the locally maximally dense arrangements of a fixed number of circles on a given flat torus, from which we determine the globally maximally dense packing(s).

The main structure that allows us to form a list of all the locally maximally dense packings for a fixed number of circles on a flat torus is the graph of a packing. Given a packing $\mathcal{P}$ on a flat torus, the \textbf{packing graph associated to }$\bm{\mathcal{P}}$, denoted $G_{\mathcal{P}}$, has geometric vertices and edges defined as follows. (This is also sometimes called a kissing or contact graph as in~\cite{musin}.) The center of each circle in the packing is associated to a vertex (with a location in the torus) of $G_{\mathcal{P}}$ and two (not necessarily distinct) vertices of  $G_{\mathcal{P}}$ are connected with an edge (with a length as a line/geodesic segment on the torus) if and only if the corresponding circles are tangent to each other.  Further, all circle-circle tangencies lead to an edge and this allows pairs of vertices to have more than one edge connecting them. Thus each packing of equal circles on a flat torus is naturally associated to an embedding of a graph on a flat torus where all the edges are equal in length. Throughout this article we allow the graphs to be multigraphs that can possibly contain multiedges and loops. 

It is important to note that it is possible to have a \textbf{free circle} or \textbf{rattler} in an optimal arrangement; that is, a circle that is not tangent to any \textit{other} circle.  For packing on a torus, a circle that is only self-tangent is still considered to be free. The work of Schaer~(\cite{schaer10}) on a unit square,  Musin~(\cite{musin}) on a square torus, and Melissen~(\cite{equilateralTriangle}) on an equilateral triangle, prove that at least one of the globally optimally dense packings of seven equal circles in each of these domains contains a free circle. In the present work we discover a family of locally maximally dense equal circle packings with four circles that admits a fifth circle that is free.  See the comments after Theorem~\ref{thm:connelly} in Section~\ref{sec:rigid} for more discussion and the top row of Figure~\ref{fig:lmdwfree} for an image.  

\subsection{Overview of the proof of Theorem~\ref{thm:MainResult}}\label{subsection:overview} 
By viewing the associated packing graph of an optimal packing as type of tensegrity framework (an embedded graph with additional structure) on a torus, Connelly (\cite{connelly1}) has proven a lower bound on the number of circle-circle tangencies.  His result implies that in an optimal packing of $n$ equal circles (without any free circles) on a flat torus there must be $2n-1$ or more edges in the associated packing graph. Using the fact that each circle in a torus can be tangent to at most six circles, we can deduce that in an optimal packing of $n$ equal circles (without any free circles) on a flat torus there must be $3n$ or fewer edges.  This means that there are a finite number of graphs that can be associated to an optimal packing of $n$ equal circles on a flat torus. By studying all of these possible graphs, and all the possible ways they can be embedded on a torus, we are led to all the locally (and globally) optimally dense packings.  

More specifically, this leads us to the following outline of the proof of Theorem~\ref{thm:MainResult}. (Note that the case of $n=2$ is handled differently in Subsection~\ref{twoequalcircles} and those packings with self-tangent circles are handled separately in Section~\ref{sec:selftangent}.)
\begin{enumerate}
\item Step one is divided into two parts and the details are provided in Section~\ref{sec:rigid}.
\begin{enumerate}
\item\label{steponea} Make a list of all the possible graphs that could be associated to an optimal packing of $n=3$ or $n=4$ circles on a flat torus. That is, list all the combinatorially distinct graphs with $n=3$ or $n=4$ vertices and between $2n-1$ and $3n$ edges. There are 862 of these combinatorial graphs. 
\item\label{steponeb} Not all of the combinatorial graphs from step~\ref{steponea} can be associated to an optimal packing. We use Proposition~\ref{prop:3tangencies} to eliminate all but 23 of the combinatorial graphs.  This is shown in Table~\ref{tab:data} and all 23 of them are shown in Figure~\ref{fig:combinatorialGraphs}.  
\end{enumerate}

\item Step two is divided into two parts and the details are provided in Section~\ref{sec:TopologicalGraphTheory}.
\begin{enumerate}
\item\label{steptwoa} For each of the combinatorial graphs remaining from step~\ref{steponeb}, we use techniques from topological graph theory (specifically Edmond's permutation technique \cite{edmonds}) to enumerate all the possible 2-cell embeddings of it onto a topological torus.  Note that a combinatorial graph can embed in multiple ways. There are 103 of these embedded combinatorial graphs (or simply embedded graphs).
 
\item\label{steptwob} Not all of these embedded combinatorial graphs can associated to an equal circle packing (optimal or not). We use Propositions~\ref{prop:forbiddenFacePatterns} and~\ref{prop:edgeForce} to eliminate all but 27 of the embedded combinatorial graphs. This is shown in Table~\ref{tab:remaining} and all 27 of them are shown in Figures~\ref{fig:toroidalEmbeddings3} and~\ref{fig:toroidalEmbeddings4}.   
\end{enumerate}

\item For each of the embedded combinatorial graphs remaining from step~\ref{steptwob} we must answer some questions: Does the embedded graph correspond to an equal circle packing? If so, is the packing optimal or not? If so, is the packing locally or globally optimal? As discussed in Section~\ref{sec:nonoptimal}, 15 of the embedded combinatorial graphs either do not correspond to an equal circle packing or, if they do, the corresponding packing is locally, but not globally optimal.  Section~\ref{sec:Existence} handles the remaining 12 embedded combinatorial graphs that are associated to the globally optimal packings. In this section, explicit coordinates for the locations of centers of the circles are given and, thus,  the optimal radii at all locations in the moduli space can be determined.
\end{enumerate}
Therefore, because we started with all possible combinatorial graphs that could be associated with an optimal packing and we embedded those combinatorial graphs in all possible ways, we have determined all the globally (and locally) packings on any torus. Note: As a torus is compact and the packing radius determined by a collection of $n$ points on a torus is continuous, we know that there is a globally optimal packing of $n$ equal circles on that torus.

\section{Results From Rigidity Theory }\label{sec:rigid}
 
We now discuss some tools from rigidity theory, including some useful results from~\cite[Section~3]{dickinson3} and~\cite[Section~3]{musin}. See these references for more details. In particular, we state propositions that establish an upper and lower bound on the number of edges a packing graph associated to an optimal packing must contain.  This allows us to create a short finite list (Figure~\ref{fig:combinatorialGraphs}) of combinatorial multigraphs without loops (loops are handled in Section~\ref{sec:selftangent}) each of which is a candidate to be the packing graph of an optimal packing.  Afterwards we complete Step 1 of the proof (following the overview of the proof found in Subsection~\ref{subsection:overview}).

\subsection{Tools From Rigidity Theory}
Rigidity theory involves the study of tensegrity frameworks.  For our purposes, we specialize to strut tensegrity frameworks, which are essentially graphs embedded in a Riemannian manifold with some additional structure.  To be a strut tensegrity framework, each edge in the graph is not allowed to decrease in length as the location of its endpoint vertices change. We can view the (embedded) packing graph of a circle packing as a strut framework on a torus.  This is appropriate because as we move the vertices of a circle packing graph to try and improve the density, we want the length of the edges to either increase (or remain unchanged), in order to possibly increase (or maintain) the density. The motions (if any) of the vertices that respect the distance constraints between the vertices are called flexes. A flex that is not induced from a family of rigid motions of the torus is called a non-trivial flex. If there are no non-trivial flexes, then the strut framework is called rigid. Combining the rigidity theory ideas with the ideas of circle packing, we have the following (for a complete proof see~\cite[Prop. 3.1]{dickinson2}). 

\begin{proposition} \label{prop:converse}
If the strut tensegrity framework associated to a circle packing $\mathcal{P}$ is rigid, then the circle packing $\mathcal{P}$ is locally maximally dense.
\end{proposition}

The flexes of a framework can be linearized and will lead to another notion of rigidity that will turn out to be equivalent in our case.  If you consider the time zero derivative of a flex at each vertex, then you obtain a collection of vectors. This collection of vectors must satisfy a system of linear homogeneous (strut) inequalities which result from the flex respecting the distance constraints. A collection of vectors that satisfy the strut inequalities and are not the time zero derivative of a family of rigid motion of the torus at each vertex is called a non-trivial infinitesimal flex.  If there are no non-trivial infinitesimal flexes, then the strut framework is called infinitesimally rigid. The connection between infinitesimal rigidity and rigidity of a framework is well studied and Connelly (\cite{connelly1}) has proven that a strut tensegrity framework is rigid if and only if it is infinitesimally rigid. Determining the infinitesimal rigidity of an arrangement is the same as the feasibility part of linear programming, is straightforward to check, and, with Proposition~\ref{prop:converse},  enables us to easily check when a given circle packing is locally maximally dense.

We are now in a position to state the main result from Connelly~(\cite{connelly1}), specialized to the context of flat tori, which is almost the converse of Proposition~\ref{prop:converse}.
\begin{theorem}[Connelly]
\label{thm:connelly}
If $\mathcal{P}$ is a packing that is locally maximally dense on a flat torus,  then there is a sub-packing $\mathcal{Q}$ of $\mathcal{P}$ such that the associated strut tensegrity framework is infinitesimally rigid and the circles not in $\mathcal{Q}$ (possibly an empty set) are prevented from increasing their radius (i.e. are free circles).
\end{theorem}
Motivated by Theorem~\ref{thm:connelly} we now discuss those locally maximally dense packing with disconnected packing graphs. 

\subsubsection{Optimal Packings With A Disconnected Packing Graph Or That Admit Another Circle}
If we remove any free circles from a locally maximally dense arrangement, then we obtain a locally maximally dense packing for fewer circles in the flat torus.  Conversely, if there is room for another circle (free or not) in a globally maximally dense packing of $n$ circles, then adding this circle gives us an arrangement realizing the globally maximal density for $n+1$ circles. Heppes~(\cite{heppes})and Melissen ~(\cite{melissen}) exploited this when they noticed that their globally optimal packings for 3 circles in the boundary of region $R1_3$ that is along the $y$ axis admit another circle. That is, the globally maximally dense packings along this edge have room for a fourth circle (creating an additional 4 tangencies so the circle is not free). This leads to the globally maximally dense arrangements in the corresponding locations in region $R1_4$ on the $y$ axis for 4 circles.

The present work gives more examples of this phenomenon. The family of packings for 3 circles that is globally optimally dense for region $R1_3$ also extends to a locally optimally dense packing on the left of the $y$ axis (this is the light gray region in the middle of Figure~\ref{fig:ModuliSpaces}).  It turns out that this locally optimal two parameter family of 3 circles packing admits a fourth circle (with 4 additional tangencies) that becomes globally optimally dense packings for 4 circles. See the top row of Figure~\ref{fig:extended}. With the addition of this circle, these packings occupy region $R1_4$ on the right of Figure~\ref{fig:ModuliSpaces}.   

The only time that a free circle (or an additional circle) might be able to be added to an arrangement is when there is a face in the packing graph consisting of seven or more edges (\cite[Prop. 3.6]{musin}). There are two cases in the present work where this happens in an optimal packing.
\begin{itemize}
\item The optimal packing occupying region $R1_3$ (where it doesn't admit the additional circle) and the adjacent gray region (seen in the middle of Figure~\ref{fig:ModuliSpaces} which is discussed in the paragraph above) for 3 circles.
\item In one locally (but not globally) maximally dense family of packings of 4 circles. This does admit a free circle. See the top row of Figure~\ref{fig:lmdwfree}. 
This packing is a locally optimally dense packing of 5 circles which is beyond the scope of the present work. 
\end{itemize}
It should be noted that there are two arrangements of 3 circles that contain a regular hexagonal face with room for an additional circle (with 6 additional tangencies). Adding the fourth circle to these arrangements yields the triangular close packing with 4 circles on tori  corresponding to the unfilled points on the vertical axes on the right side of Figure~\ref{fig:ModuliSpaces}.  There are a handful of locally optimal 4 circle packings 
 with a regular hexagon and adding a circle leads to the triangular close packing on 5 circles. 

In this article, we will determine all the locally maximally dense arrangements for 3 and 4 equal circles without free circles.  This, coupled with the discussion here, means that we will have created an exhaustive list of locally maximally dense packings.  Therefore, for the remainder of this article, we assume that all of our graphs are connected. As noted earlier, Przeworski~\cite[Theorem 2.3]{przeworski} has solved this problem for 2 equal circles and, for completeness we analyze this case using our methods (see Subsection~\ref{twoequalcircles}) and we reach the same conclusion as Przeworski.

\subsubsection{Characterizing The Packing Graphs Associated To Optimal Packings}
Now we observe that we can find a lower bound on the number of edges (and their arrangement) incident to a vertex in the packing graph associated to a locally maximally dense packing with no free circles.

\begin{proposition}
\label{prop:VertexLowerBound}
If $\mathcal{P}$ is a locally maximally dense packing of circles on a flat torus with no free circles, then no circle in $\mathcal{P}$ has its points of tangency contained in a closed semi-circle. In particular, every circle is tangent to at least three circles.
\end{proposition}

\begin{proof}
If there were such a circle in a locally maximally dense packing, then the packing graph would not be infinitesimally rigid violating Theorem~\ref{thm:connelly}. For more details see~\cite[Section~3]{dickinson3}.
\end{proof}

Connelly~(\cite{connelly2}) proves a lower bound on the number of edges that a packing graph must contain in order for the associated packing to be locally maximally dense.  

\begin{proposition}[Connelly]
\label{prop:minimum}
If $\mathcal{P}$ is a locally maximally dense packing of $n$ circles on a flat torus with no free circles,  then the packing graph associated to $\mathcal{P}$ contains at least $2n-1$ edges.
\end{proposition}

In the case of 3 and 4 equal circles (unlike the case of 5 or more circles on a square torus, see~\cite[Prop. 4.4]{dickinson2}) it is possible that a circle can be tangent to another circle in two different ways. However, we can eliminate the possibility of two circles being tangent in three (or more) ways in the case of 4 (or more) circles and restrict the arrangements in the case of 3 circles to the triangular close packing. 

\begin{figure}[h]
\includegraphics{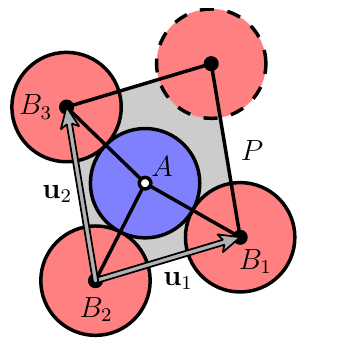} 
\caption{A packing in which circle $A$ is tangent to circle $B$ in at least three ways.}\label{fig:Eliminate3MoreTangency}
\end{figure}
\newpage

\begin{proposition}\label{prop:3tangencies}
Given an equal circle packing on a flat torus with $n$ circles.
\begin{itemize}
\item If $n \geq 4$, then no two circles can share three or more tangencies. 
\item If $n=3$, then no two circles can share four or more tangencies. Further if there are two circles that share three tangencies, then every circle is tangent to six circles (i.e. the packing lifts to the triangular close packing).
\end{itemize}
\end{proposition}
\begin{proof}
Suppose circle $A$ is tangent to circle $B$ in 3 (or more ways) on a torus that is the quotient of the Euclidean plane by lattice $\Lambda$. Lift circle $A$ into the Euclidean plane and consider the triangle formed by the centers of the lifts of circle $B$ that are tangent to the lift of $A$. Call these centers $B_1$, $B_2$, and $B_3$. 
Let $\textbf{u}_1$ ($\textbf{u}_2$) be the vector of $\Lambda$ that connects $B_2$ and $B_1$ ($B_2$ and $B_3$) respectively. Let $P$ be the $\Lambda$-lattice parallelogram formed by  $B_1$, $B_2$, $B_3$ and $B_2 + \textbf{u}_1 + \textbf{u}_2$. See Figure~\ref{fig:Eliminate3MoreTangency}. 
By Pick's theorem the area of $P$, $A_P$, is equal to $A_{FD}(i+\frac{b}{2}-1)$ where $A_{FD}$ is the area of a fundamental domain of lattice $\Lambda$, $i$ ($b$) is the number of lattice points interior (on the boundary) of $P$. Using  the minimum number of lattice points of $\Lambda$ inside and on the boundary of $P$, we have that $A_P \geq A_{FD}$.  Further if $d$ is the common diameter of the circles, then the geometry of the packing implies that the maximum area of $A_P$ is $\frac{3\sqrt{3}}{2}d^2$
and
we have the following bounds on the density, $\rho$, of the packing.
\begin{equation*}
\frac{\pi}{\sqrt{12}} \geq \rho = \frac{\text{Area covered by $n$ circles}}{A_{FD}} \geq \frac{n \pi \left(\frac{d}{2}\right)^2}{A_P} \geq \frac{n \pi \left(\frac{d}{2}\right)^2}{\frac{3\sqrt{3}}{2}d^2} \geq \frac{n\pi}{3\sqrt{12}}
\end{equation*}
For $n\geq4$ this is a contradiction and for $n=3$ all of the inequalities become equalities and the density of the packing implies that the arrangement must lift to the triangular close packing.  This also implies that it is impossible for a packing of three circles on a flat torus to have a pair of circles tangent in four or more ways. 
\end{proof}

These three propositions and the observation that in the triangular close packing (the most dense packing of equal circles in the plane) each circle is tangent to six others, helps us list important properties of the packing graphs of locally maximally dense arrangements $n$ circles on the torus. This is summarized in the following proposition.

\begin{proposition} \label{prop:GraphProperties} Given a locally maximally dense packing, $\mathcal{P}$, of $n \geq 3$ equal circles without any free or self-tangent circles on a flat torus, the packing graph $G_\mathcal{P}$ satisfies the following three conditions.
\begin{enumerate}
\item\label{blah} $G_\mathcal{P}$ is connected, contains no loops, and contains at least $2n-1$ and at most $3n$ edges.
\item Every vertex of $G_\mathcal{P}$ is connected to at least three and at most six others.
\item No pair of vertices of $G_\mathcal{P}$ is connected by 3 or more edges, except if $n=3$ and the packing is the triangular close packing implying that all three vertices of $G_\mathcal{P}$ have degree 6.
\end{enumerate}
\end{proposition}

\subsection{Step~1 of the proof of~\ref{thm:MainResult}}
\begin{figure}
 \begin{tabular}{ccccc}
 &\includegraphics{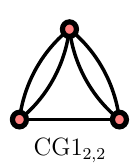} & \includegraphics{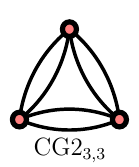} & \includegraphics{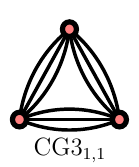} &\\
 %

%
  \includegraphics{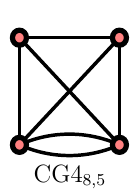} & \includegraphics{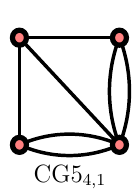} &  \includegraphics{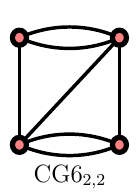} &\includegraphics{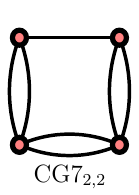} & \includegraphics{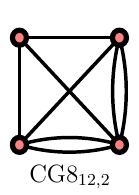}    \\
 \includegraphics{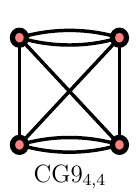} & \includegraphics{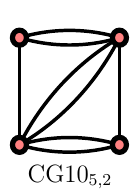} &\includegraphics{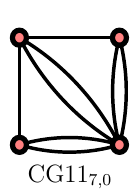}  & \includegraphics{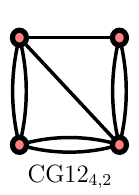} &\includegraphics{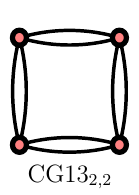}  \\
  \includegraphics{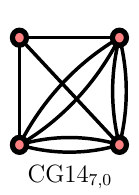} & \includegraphics{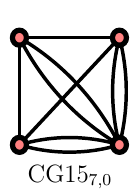} & \includegraphics{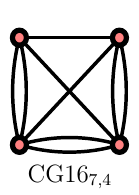} &  \includegraphics{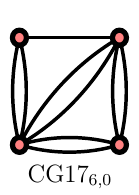} &\includegraphics{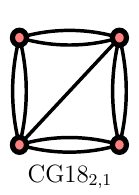} \\
 \includegraphics{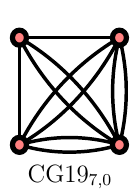} &\includegraphics{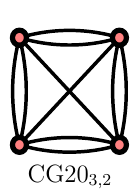} &  \includegraphics{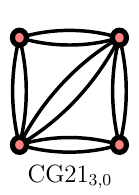} & \includegraphics{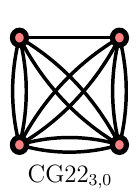}&\includegraphics{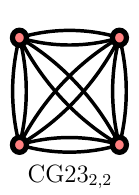} \\
 \end{tabular}
 \caption{The 23 combinatorial graphs (CG) with three or four vertices satisfying Proposition~\ref{prop:GraphProperties}. The first number in the subscript is the number of distinct unlabeled and unoriented toroidal embeddings. The second number in the subscript is the number of distinct unlabeled and unoriented toroidal embeddings remaining after applying Proposition~\ref{prop:forbiddenFacePatterns}.  (Applying Prop.~\ref{prop:edgeForce} eliminates more.) The remaining toroidal embeddings are shown in Figures~\ref{fig:toroidalEmbeddings3} and~\ref{fig:toroidalEmbeddings4}. See Section~\ref{sec:TopologicalGraphTheory} for more details about these subscript numbers.}  \label{fig:combinatorialGraphs}
 \end{figure}

The number of multigraphs (we allow multiple edges between vertices because these correspond to pairs of circles tangent in multiple ways) with a fixed number of vertices and edges is well studied.  Using Dr. Gordon Royle's data posted on the web~(\cite{gordon}) or the program Nauty (\cite{nauty} with the gtools \texttt{geng} and \texttt{multig}) we can make a list of the combinatorial graphs that satisfy the first condition of Proposition~\ref{prop:GraphProperties}. (This is step~\ref{steponea} of the overview of the proof found in Subsection~\ref{subsection:overview}.) Both the posted data and the results of Nauty yielded the same number of graphs shown in the first column in Table~\ref{tab:data}.  Using Wolfram's Mathematica (\cite{Mathematica}) we wrote a routine to remove the combinatorial graphs that didn't satisfy conditions two and three of Proposition~\ref{prop:GraphProperties}.  See Figure~\ref{fig:combinatorialGraphs} for visualizations of the remaining combinatorial graphs. (This is step~\ref{steponeb} of the overview of the proof found in Subsection~\ref{subsection:overview}.) 

\begin{table}
\renewcommand{\arraystretch}{2.5}
\begin{tabular}{|c|c|c|c|c|c|} \hline
                 & \parbox[c]{1.25in}{Number of graphs satisfying condition 1 of Prop.~\ref{prop:GraphProperties}} & \parbox[c]{1.25in}{Number of graphs satisfying conditions 1 \& 2 of Prop.~\ref{prop:GraphProperties}}& \parbox[c]{1.25in}{Number of graphs satisfying all conditions of Prop.~\ref{prop:GraphProperties}} \\ \hline
3 Vertices & 37 & 10 & 3 \\ \hline
4 Vertices & 825 & 102 &20 \\ \hline
\end{tabular}
\vspace{.1in}
\caption{The number of combinatorial graphs remaining after parts of Proposition~\ref{prop:GraphProperties}. See Figure~\ref{fig:combinatorialGraphs} for visualizations. }\label{tab:data}
\end{table}

 \section{Packings with self-tangent circles}\label{sec:selftangent}
\begin{figure}
\includegraphics{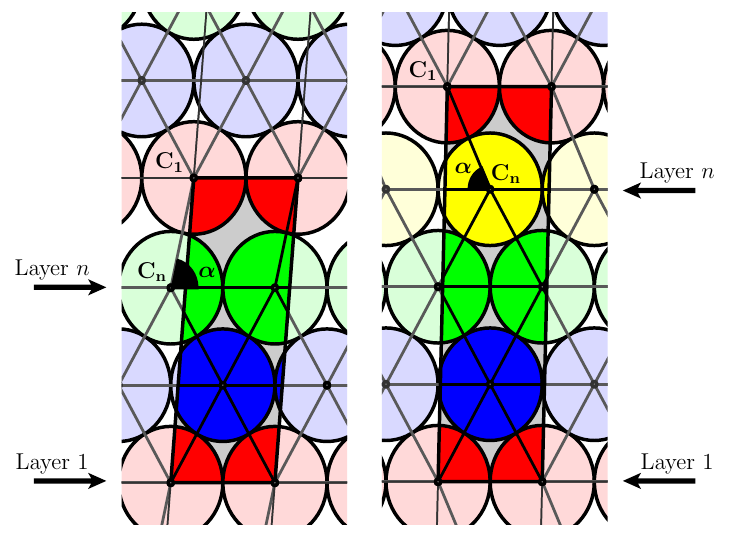}
\caption{The lattice associated to an optimal packing with an odd number (left) or an even number (right) of equal circles with packing radius $\frac{1}{2}$ with no free circles. Notice how all the circles are self-tangent and each forms a layer so that towards the bottom of the torus the layers are arranged to form part of a triangular close packing.}\label{fig:3}
\end{figure}

In this section, we completely characterize all optimal packings of 3 or 4 equal circles whose packing graph contains a loop. This leads us to justify the assumption used outside of this section that all packing graphs do not contain a loop.  On the torus, locally maximally dense packings with free circles are closely related to packings that have self-tangent circles. This is because if you have a torus with basis $\textbf{v}_1 = \langle 1, 0\rangle$ and $ \textbf{v}_2 =\langle x, y\rangle$ where $x^2+y^2\geq 1$, $y>0$, and $0 \leq x \leq \frac{1}{2}$ and the length of $\textbf{v}_2$ is long enough, all the circles can achieve the maximum radius of $\frac{1}{2}$ and become self-tangent and free. In this section we determine the boundary between where all the circles are self-tangent (see Figure~\ref{fig:3}) and where none of the circles are self-tangent.  We show that the region allowing self-tangent circles is the only one where packing graphs with a loop are realizable; this means that we have characterized all packings with a loop in them and need not consider combinatorial graphs with loops outside of this section.  Throughout this section we assume that the torus is the quotient of the Euclidean plane by the vectors $\textbf{v}_1$ and $\textbf{v}_2$ that satisfy the restrictions given in this paragraph.  We begin with some observations about packings that contain self-tangent circles.

\begin{proposition} \label{prop:1}
Let $\mathcal{P}$ be a packing of equal circles on a flat torus. The following four statements are equivalent.
\begin{enumerate}
\item There exists a circle in $\mathcal{P}$ that is self-tangent.
\item The combinatorial multigraph associated to $\mathcal{P}$ contains a loop. 
\item The common radius of the circles in $\mathcal{P}$ is $\frac{1}{2}$.
\item All the circles of $\mathcal{P}$ are self-tangent.
\end{enumerate} 
\end{proposition}

\begin{proof}  All statements follow from the fact that this is an equal circle packing and the observation that $\textbf{v}_1$ is a shortest vector of length one so lifts of circles that differ by this lattice element must have radius $\frac{1}{2}$. In this case all circles in the packing must be self-tangent.
\end{proof}

Observe that Proposition~\ref{prop:1} implies that if there is a loop in the combinatorial graph associated to a packing $\mathcal{P}$ (or the packing radius is $\frac{1}{2}$), then $\mathcal{P}$ is locally maximally dense. This is because $\frac{1}{2}$ is the absolute upper bound of the radius in any packing on the tori we are considering so any nearby packing has density less or equal to that of $\mathcal{P}$. 

Now we examine the structure of packings with a self-tangent circle and consider the cases where such a packing contains or does not contain a free circle. The intuition behind the first part of Proposition~\ref{prop:2} is to first note that all the circles in the packing have radius $\frac{1}{2}$ and therefore each circle is self-tangent, so it and its lifts create a ``layer" in the flat torus -- see Figure~\ref{fig:3}. Next, in order to fit the packing on the ``smallest" flat torus, the layers must be stacked so that most circles are tangent to the circle above and below in two different ways forming a section of triangular close packing at the ``bottom" of the torus. 

\begin{proposition} \label{prop:2}
Let $\mathcal{P}$ be a packing of $n  \geq 2$ equal circles on a flat torus where the combinatorial multigraph associated to $\mathcal{P}$ contains a loop. 
\begin{enumerate}
\item \label{one} If there are no free circles in $\mathcal{P}$, then there exists $\alpha$ with $\frac{\pi}{3} \leq \alpha \leq \frac{\pi}{2}$ such that $\textbf{\emph{v}}_2 = \langle x, \frac{n-1}{2}\sqrt{3}+\sin{(\alpha)} \rangle$ where $x=\frac{1}{2}-\cos{(\alpha)}$ for $n$ even or $x= \cos{(\alpha)}$ for $n$ odd.

\item \label{two} If there is a free circle in $\mathcal{P}$,  then all circles in $\mathcal{P}$ are free and there exists $\frac{\pi}{3} \leq \alpha \leq \frac{\pi}{2}$ such that $\textbf{\emph{v}}_2 = \langle x, y\rangle$ where $x=\frac{1}{2}-\cos{(\alpha)}$ (for $n$ even) or $x=\cos{(\alpha)}$ (for $n$ odd) and $y > \frac{n-1}{2}\sqrt{3}+\sin{(\alpha)}$.

\end{enumerate} 
\end{proposition}

\begin{proof}
First we begin by showing that there exists a cyclic ordering of the circles $C_1, C_2, C_3, \ldots, C_n$ where circle $C_i$ is tangent to itself and can only be tangent to $C_{i-1}$ or $C_{i+1}$ for $1 \leq i \leq n$ (indices of the circles counted modulo $n$).
To see this notice that as the combinatorial multigraph associated to $\mathcal{P}$ contains a loop, by Proposition~\ref{prop:1}, all circles are self-tangent. For each circle in the packing, this means that the toroidally embedded packing graph contains an essential cycle (a closed path not homotopic to a point) of length one consisting only of the edge corresponding to the self-tangency. Cutting the torus along these essential cycles results in $n$ cylinders. This establishes a cyclic ordering of the circles.  As all tangencies in this packing must occur in some cylinder, we know that, beside being tangent to itself, $C_i$ can only be tangent to the circles forming the edges of the cylinders with the essential cycle associated to $C_i$, namely $C_{i-1}$ and $C_{i+1}$. 

To establish item (\ref{one}) of the proposition, assume that there are no free circles in $\mathcal{P}$ and by the remark immediately preceding this proposition, $\mathcal{P}$ is locally maximally dense.  Theorem~\ref{thm:connelly} with the assumption that there are no free circles implies that the entire strut tensegrity framework associated $\mathcal{P}$ must be infinitesimally rigid (not just some subgraph). This implies that each circle in the packing must have at least three tangencies not including the self-tangency, because the strut inequality from the self-tangency is trivial. The only way for this to be true is if each circle $C_i$ is tangent at least once (and at most twice) to each of $C_{i-1}$ and $C_{i+1}$.  

Now we observe that there can be at most one $i$ such that $C_i$ is tangent once to $C_{i+1}$.  Suppose not, then there exists two natural numbers, $i$ and $j$ ($i < j \leq n$), such that $C_i$ is tangent only once to $C_{i+1}$, $C_j$ is tangent only once to $C_{j+1}$ and for each natural number $m$, $i+1 \leq m \leq j-1$, $C_m$ is tangent twice to $C_{m+1}$. In this case there is a non-trivial infinitesimal flex  (i.e. a non-trivial assignment of vectors to the circle centers that satisfies the strut inequalities -- See Section~\ref{sec:rigid} for more details) that is a non-zero constant vector at the centers of $C_m$ ($i+1 \leq m \leq j-1$) and the zero vector at the other circle centers. For the non-zero constant, choose a vector that makes an obtuse angle with both the single edge between $C_i$ and $C_{i+1}$ and the edge between $C_j$ and $C_{j+1}$. Such a choice always exists except in the case when these two edges are parallel (in which case you would choose the non-zero vector to be perpendicular to both). Hence there is at most one circle $C_i$ that is tangent once to $C_{i+1}$.

For the purpose of determining the lattice that defines the torus assume, by renumbering if necessary, that the possibly single tangency between circles occurs between $C_n$ and $C_1$. By lifting the packing to the plane and by reflecting if necessary we can assume that the angle, $\alpha$, between the horizontal and the edge between $C_n$ and $C_1$ (measured clockwise for $n$ even and measured counterclockwise for $n$ odd) is between $\frac{\pi}{3}$ and $\frac{\pi}{2}$. In this case $\textbf{v}_2$ has the form given in item~(\ref{one}) of this proposition. See Figure~\ref{fig:3}.  

For item~(\ref{two}), suppose that $\mathcal{P}$ contains a free circle, $C_i$.  As this circle is free we can translate it to break all tangencies (if any) with this circle (except the self-tangency).  Next, we can translate (perpendicular to $\textbf{v}_1$) all other circles a little bit using the space created when all tangencies with $C_i$ were broken. In this way we see that all circles are free in this packing.  To determine the lattice of the torus for this packing notice that as the circles are free we can move them into the positions so that there are two tangencies between $C_i$ and $C_{i+1}$ for $1 \leq i < n$.  This results in a situation similar to that which is pictured in Figure~\ref{fig:3} except that there are no tangencies between circle $C_n$ and $C_1$. Notice that the angle $\alpha$ is the angle such that the $x$ coordinate of $\textbf{v}_2$ is either $\cos(\alpha)$ ($n$ odd) or $\frac{1}{2}-\cos{(\alpha)}$ ($n$ even). The bound on $y$ coordinate follows.  
\end{proof}

Finally we use the constructions, results, and observations that arise in the rest of this article and Proposition~\ref{prop:2} to observe that all locally and globally optimal packings on tori in the region described in item~(\ref{two}) must have a loop in the packing graph. Therefore we have completely characterized the tori on which there exists a locally maximally dense packing that contains a self-tangent circle for 3 or 4 equal circles.  Notice that the descriptions of the regions for odd and even numbers of circles given in Propostion~\ref{prop:3} agree with the descriptions of regions $R2_2$, $R3_3$, and $R4_4$ given in Table~\ref{tab:regions}. 

\begin{proposition} \label{prop:3}
Let $\mathcal{P}$ be a locally maximally dense packing of $n=3$ or $n=4$ equal circles on a flat torus. If the flat torus is the quotient of the plane with respect to a lattice with basis vectors $\text{v}_1 = \langle 1,0 \rangle$ and $\textbf{\emph{v}}_2 = \langle x, y\rangle$ where there exists $\frac{\pi}{3} \leq \alpha \leq \frac{\pi}{2}$ such that $x=\frac{1}{2}-\cos{(\alpha)}$ (for $n$ even) or $x=\cos{(\alpha)}$ (for $n$ odd) and $y > \frac{n-1}{2}\sqrt{3}+\sin{(\alpha)}$, then there is loop in the combinatorial multigraph associated to $\mathcal{P}$ and $\mathcal{P}$ contains only free circles.
\end{proposition}
\begin{proof}
Suppose not, then there is no loop in the associated combinatorial multigraph. For $n =3,4$ all possible multigraphs without loops were considered as packing graphs in the remainder of this article and none of them embed locally maximally densely on a torus with basis $\textbf{v}_1$ and $\textbf{v}_2$ as described in this proposition. Therefore there must be a loop in the combinatorial multigraph. If $\mathcal{P}$ contains no free circles, then Proposition~\ref{prop:2} item~(\ref{one}) gives the form of the basis vectors for the torus, but these are incompatible with the basis vectors for the torus described in this proposition, therefore there must be at least one free circle. Proposition~\ref{prop:2} item~(\ref{two}) implies that all the circles in $\mathcal{P}$ are free. 
\end{proof}
Propositions~\ref{prop:3} and~\ref{prop:2} implies that we do not have to consider loops in the packing graph explorations in the remaining sections of this article.

\begin{figure}
\includegraphics{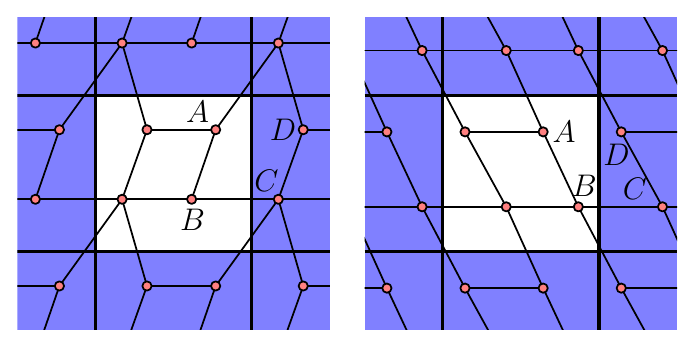}
\caption{Toroidal embeddings of CG5 (left) and CG4 (right) that do not correspond to any equal circle packing graph on any torus. The edges $\overline{AB}$ and $\overline{CD}$ are parallel because of the rhombus in the embedded graph, but there is no edge $\overline{AD}$, so neither of these can be packing graphs associated to an equal circle packing.}\label{fig:eliminategraphs}
\end{figure}

\section{Results From Topological Graph Theory}\label{sec:TopologicalGraphTheory}

This section is broken into two parts mirroring two parts of step 2 of the overview of the proof from Subsection~\ref{subsection:overview}. First, we apply previously known techniques for determining all the 2-cell embeddings of a combinatorial multigraph onto a topological torus. Second, we apply several propositions that prohibit a toroidally embedded graph from being associated with an optimal equal circle packing. This leaves us with a short list of toroidally embedded graphs that potentially could be associated with an optimal equal circle packing. See Figures~\ref{fig:toroidalEmbeddings3} and~\ref{fig:toroidalEmbeddings4}.

\subsection{Step~\ref{steptwoa} of the proof of Theorem~\ref{thm:MainResult}}
A graph embedded in a torus is a 2-cell embedding if each connected region (face) determined by removing the embedded graph from the torus is homeomorphic to an open disk. As noted in~\cite{dickinson3,musin} if a packing on the torus is locally maximally dense, then the associated packing graph is a 2-cell embedding. One tool for enumerating all the 2-cell embeddings of a graph on any surfaces is Edmonds' permutation technique~(\cite{edmonds}) which is outlined in~\cite[Section~5]{dickinson3}. (We use the same C++ implementation, adapted for multigraphs, in this article.) Essentially this is a brute force technique where, at each vertex, all possible orderings of the adjacent vertices are considered (call rotation schemes). Once a face walking algorithm is executed and the Euler characteristic is computed, all the possible toroidal embeddings can be selected. As the graphs we are dealing with are small, this brute force method converges. There are more efficient algorithms, see~\cite[Sect. 13.4]{kocay3}. The software written by Kocay~(\cite{kocay1}) was instrumental in visualizing and computing the these embedded graphs.  The paper, \textit{Embeddings of small graphs on the torus} by Gagarin \textit{et al.} ~(\cite{kocay4}), present an independently created table listing the number of toroidal embeddings of all vertex-transitive graphs (with fewer than 12 vertices) onto the torus. The results of our program were checked against this table and, when our brute force program converged, it agreed with the results in this paper.

\subsection{Step~\ref{steptwob} of the proof of Theorem~\ref{thm:MainResult}}
Once we have a list of all possible toroidal embeddings we ask if the embedded graph could be associated to an optimal equal circle packing. The following two propositions eliminate many of these potential packing graphs. 

\begin{proposition}\label{prop:forbiddenFacePatterns}
If a graph embedded on a torus contains a vertex surrounded by any one of the following face patterns, then the embedded graph cannot be the graph associated to a locally maximally dense equal circle packing. The forbidden face patterns are
(1) two triangles and a polygon, 
(2) three triangles and a polygon,
(3) five triangles,
(4) four triangles and a quadrilateral,
(5) six polygons with at least one non-triangle,
(6) a triangle, a quadrilateral and a polygon,
(7) two triangles and two quadrilaterals,
(8) three quadrilaterals, or
(9) seven (or more) polygons.
\end{proposition}

For a proof of Proposition~\ref{prop:forbiddenFacePatterns} see~\cite[Prop.~6.1]{dickinson3}. Proposition~\ref{prop:edgeForce} is another useful tool for eliminating embedded graphs from being associated with any equal circle packing. It is illustrated in Figure~\ref{fig:eliminategraphs}. Table~\ref{tab:remaining} shows how these two propositions cut down on the number of embedded graphs that might correspond to an optimal equal circle packing.

\begin{proposition} \label{prop:edgeForce}
Suppose that an embedded graph corresponds to an equal circle packing and has a pair of edges, $\overline{AB}$ and $\overline{CD}$, that when lifted to the plane determine lines that are parallel. If $A$ and $D$ are on the same side of line $\overleftrightarrow{BC}$ (in the plane) and $\overline{BC}$ is another edge of the graph, then $\overline{AD}$ is also an edge in the graph.
\end{proposition}
\begin{proof}
Suppose that an embedded graph corresponds to the packing graph of an equal circle packing on a torus and that $A$, $B$, $C$, and $D$ satisfy the conditions stated in the proposition. See either part of Figure~\ref{fig:eliminategraphs}. This implies that $\overline{AB}$, $\overline{BC}$, and $\overline{CD}$ is a chain of edges each of which has length to equal the common diameter, $d$, of the circles. The geometry of this situation forces vertices $A$ and $D$ to be distance $d$ apart. This implies that the equal circles centered at these vertices are tangent and that there is an edge between them.
\end{proof}

\begin{table}
\renewcommand{\arraystretch}{2.5}
\begin{tabular}{|c|c|c|c|} \hline
&\parbox[c]{1.25in}{Number of distinct toroidal embeddings}& \parbox[c]{1.25in}{Number of embeddings remaining after Prop.~\ref{prop:forbiddenFacePatterns}} & \parbox[c]{1.5in}{Number of embeddings remaining after using Props.~\ref{prop:edgeForce} and~\ref{prop:forbiddenFacePatterns}} \\ \hline
3 Vertices & 6 & 6  & 6 \\ \hline
4 Vertices & 97 & 31 & 21 \\ \hline 
\end{tabular}
\vspace{.1in}
\caption{The number of embedded graphs that might correspond to an optimal packing after applying Propositions~\ref{prop:edgeForce} and~\ref{prop:forbiddenFacePatterns}. See Figures~\ref{fig:toroidalEmbeddings3} and~\ref{fig:toroidalEmbeddings4} for visualizations of the remaining embedded graphs.}\label{tab:remaining}
\end{table}

\begin{figure}
\begin{tabular}{ccccccc}
\includegraphics{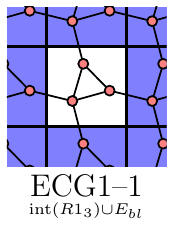} &
\includegraphics{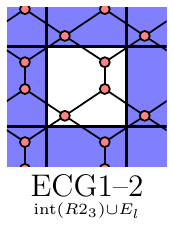} &
\includegraphics{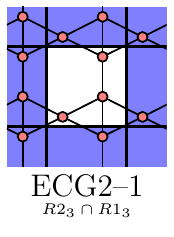} &
\includegraphics{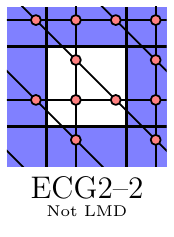} &
\includegraphics{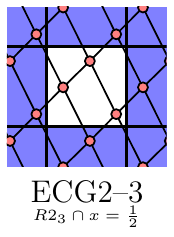}  &
\includegraphics{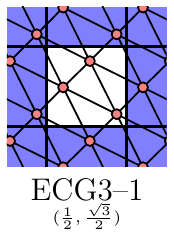} 
\end{tabular}
\caption{The 6 embedded combinatorial graphs (ECG) on three vertices that remain after applying Propositions~\ref{prop:forbiddenFacePatterns} and~\ref{prop:edgeForce}. On the first line, the number before the dash refers to the combinatorial graph (CG) that led to the embedding and the number after the dash is the embedding number.  The second line indicates more about the embedding that is useful in understanding the regions $Ri_3$ of Theorem~\ref{thm:MainResult}. If there is a globally optimal packing whose associated packing graph is the one under consideration, then the tori containing that packing occupies a region in the moduli space and it is loosely indicated. For example, ECG1--1 indicates ``int($R1_3$)$\cup E_{bl}$" which means the interior of region $R1_3$ and the bottom ($b$) and left ($l$) edges of that region. If the embedding corresponds to an equal circle packing, but it is not locally maximally dense (LMD), then it reads ``Not LMD". } \label{fig:toroidalEmbeddings3}
\end{figure}

\begin{figure}
\begin{tabular}{ccccccc}
\includegraphics{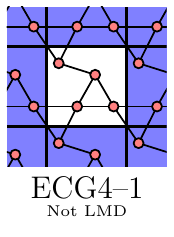} &
\includegraphics{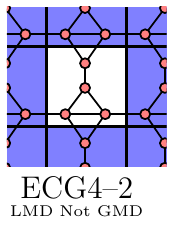} &
\includegraphics{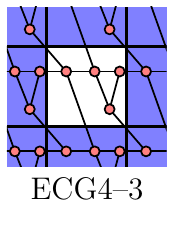} &
\includegraphics{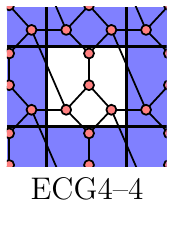} &
\includegraphics{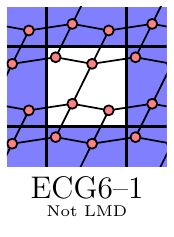} &
\includegraphics{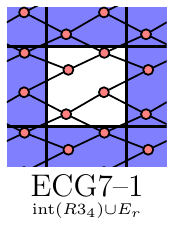}  \\
\includegraphics{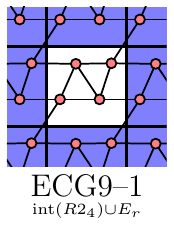} &
\includegraphics{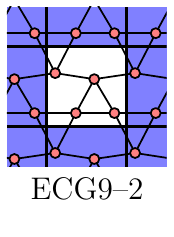} &
\includegraphics{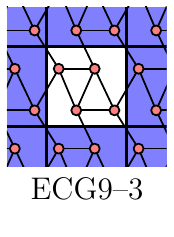} &
\includegraphics{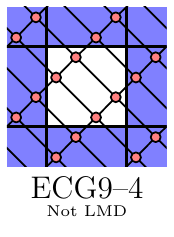} &
\includegraphics{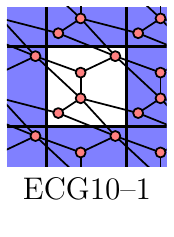} &
\includegraphics{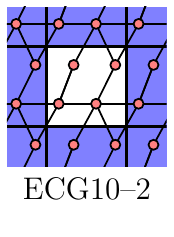} \\
\includegraphics{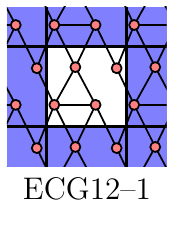} &
\includegraphics{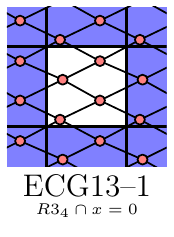} &
\includegraphics{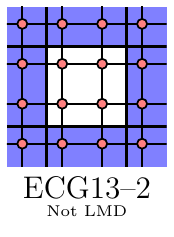} &
\includegraphics{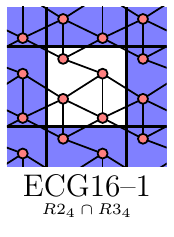} &
\includegraphics{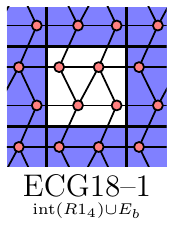} &
\includegraphics{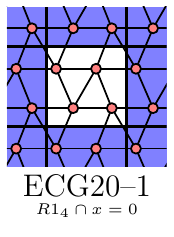} \\
\includegraphics{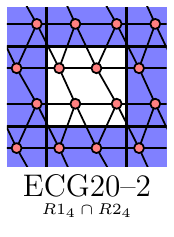} &
\includegraphics{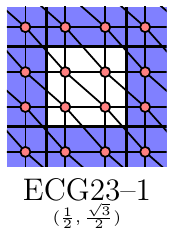} &
\includegraphics{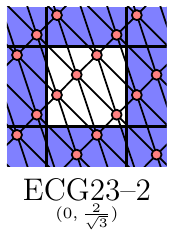} &&
\end{tabular}
\caption{The 21 embedded combinatorial graphs (ECG) on four vertices that remain after applying Propositions~\ref{prop:forbiddenFacePatterns} and~\ref{prop:edgeForce}. On the first line, the number before the dash refers to the combinatorial graph (CG) that led to the embedding and the number after the dash is the embedding number.   The second line indicates more about the embedding that is useful in understanding the regions $Ri_4$ of Theorem~\ref{thm:MainResult}. If there is a globally optimal packing whose associated packing graph is the one under consideration, then the tori containing that packing occupies a region in the moduli space and it is loosely indicated. For example, ECG23--1 indicates \protect\scalebox{.8}{$\textstyle \left(\frac{1}{2},\frac{\sqrt{3}}{2}\right)$} which means that the packing is only possible at this location in the moduli space. If the embedding corresponds to an equal circle packing, but it is not locally/globally maximally dense (LMD/GMD), then it reads ``Not LMD/GMD". If the second line is blank, then the embedding doesn't correspond to an equal circle packing on any torus.} \label{fig:toroidalEmbeddings4}
\end{figure}

\section{Non-optimal and Non-globally Maximally Dense Packings} \label{sec:nonoptimal}
In this section we discuss the details associated with first part of Step 3 of the overview of the proof found in Subsection~\ref{subsection:overview}.  We discover those embedded combinatorial graphs on three and four vertices that do not always embed in a globally maximally dense way. In summary, of the 27 embedded combinatorial graphs shown in Figures~\ref{fig:toroidalEmbeddings3} and ~\ref{fig:toroidalEmbeddings4}:
\begin{itemize}
\item Seven of the embedded graphs cannot be associated to an equal circle packing on any torus (ECG4--3, ECG4--4, ECG9--2, ECG9--3, ECG10--1, ECG10--2, ECG12--1); 
\item Five of the embedded graphs lead to a family of equal circle packing(s) on some tori, but the packings are never locally maximally dense (ECG2--2, ECG4--1, ECG6--1, ECG9--4, ECG13--2); 
\item One of the embedded graphs leads to a family of equal circle packings that are only locally maximally dense and not ever globally maximally dense (ECG4--2); and
\item Two of the embedded graphs lead to families of equal circle packings that are locally (and not globally) maximally dense on some tori and are globally maximally dense in other tori (ECG1--1, ECG9--1).
\end{itemize}

 \begin{figure}
 \includegraphics{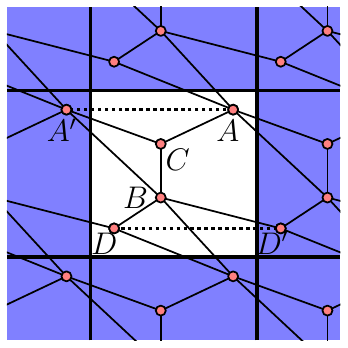}
 \caption{The embedding ECG10--1 doesn't correspond to any equal circle packing on any torus.} \label{fig:eliminateembedding}
 \end{figure}
 
Several of the embedded combinatorial graphs are not associated to an equal circle packing on any torus. The methods used to prove this are \textit{ad hoc} but usually involve showing that if all the edges in the embedding have equal length, then there are a pair of unconnected vertices that are forced to be too close together.  For example, we can eliminate the embedding 
 ECG10--1 
 in Figure~\ref{fig:eliminateembedding} from being associated to any equal circle packing on any torus by making the following observations. (Note: in the discussion below, when we say that two edges in the torus are parallel, we mean that when lifted to the Euclidean plane, they determine two lines that are parallel.) If this was associated to some equal circle packing, then all edge lengths would be equal and dashed segments $\overline{AA'}$ and $\overline{DD'}$ would be equal in length and parallel (when lifted to the plane) because their endpoints differ by the same lattice vector.  This makes triangles $\triangle AA'C$ and $\triangle DD'B$ congruent which implies that  $\overleftrightarrow{D'B}$ is parallel to $\overleftrightarrow{AC}$ or $\overleftrightarrow{A'C}$ depending on which side of $\overleftrightarrow{DD'}$ that the point $B$ is located.  If $\overleftrightarrow{D'B}$ is parallel to $\overleftrightarrow{AC}$, then Proposition~\ref{prop:edgeForce} applies to the chain of edges $\overline{AC}$, $\overline{CB}$, and $\overline{BD'}$. If $\overleftrightarrow{D'B}$ is parallel to $\overleftrightarrow{A'C}$, then Proposition~\ref{prop:edgeForce} applies to the chain of edges $\overline{CA'}$, $\overline{A'B}$, and $\overline{BD'}$. This eliminates this embedding. Embeddings  
 ECG9--2, 
 ECG9--3, 
 ECG10--2, and  
 ECG12--1  
 can be eliminated in a similar way. The arguments to eliminate the embedded combinatorial graphs 
 ECG4--3 
 and ECG4--4 
 are more involved. 
 
 \begin{figure}
 \includegraphics{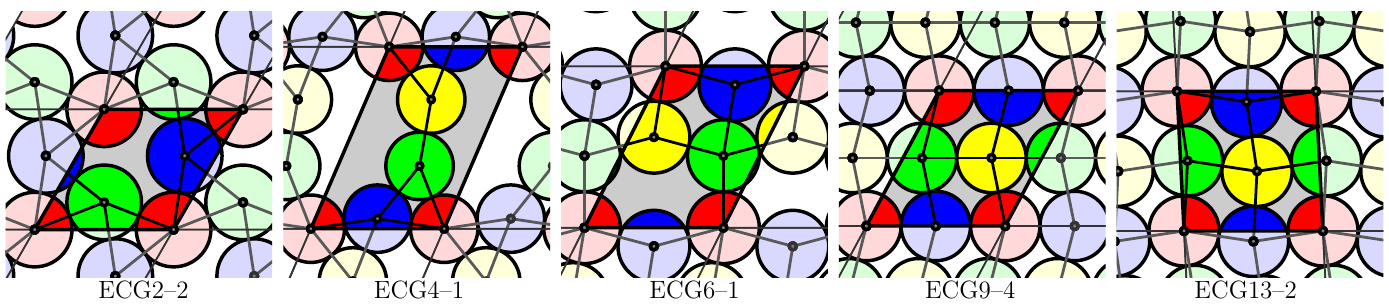}
 \caption{These are the five embedded combinatorial graphs that correspond to equal circle packings that are not locally maximally dense. Note that the two rightmost embeddings are different because the underlying combinatorial graphs are different.} \label{fig:nonstressed}
 \end{figure}
 
The embedded combinatorial graphs ECG2--2, ECG4--1, ECG6--1, ECG9--4, and ECG13--2 
correspond to circle packings that are not locally maximally dense.  To see that they correspond to equal circle packings see the circle packings in Figure~\ref{fig:nonstressed} which shows the circle packings associated to these embedded graphs. None of these packings are locally maximally dense because there is a non-trivial infinitesimal flex for each. See Section~\ref{sec:rigid} for a few details and \cite[Sec. 3]{dickinson2} for complete details and references. Roughly stated the non-trivial flex for the packing on the far left of Figure~\ref{fig:nonstressed} involves rotating the circles about the center of the equilateral triangle in the packing graph. For the remaining ones, the non-trivial flex involves fixing a circle or `row' of circles and sliding another row.  For example in the packing on the second to the left, if you fix the circle at the lower left corner of a fundamental domain, then the remaining circles can all be `slid' upwards. In the language of the article~\cite{connelly2}, these arrangements lack a proper stress.

\begin{figure}
\includegraphics{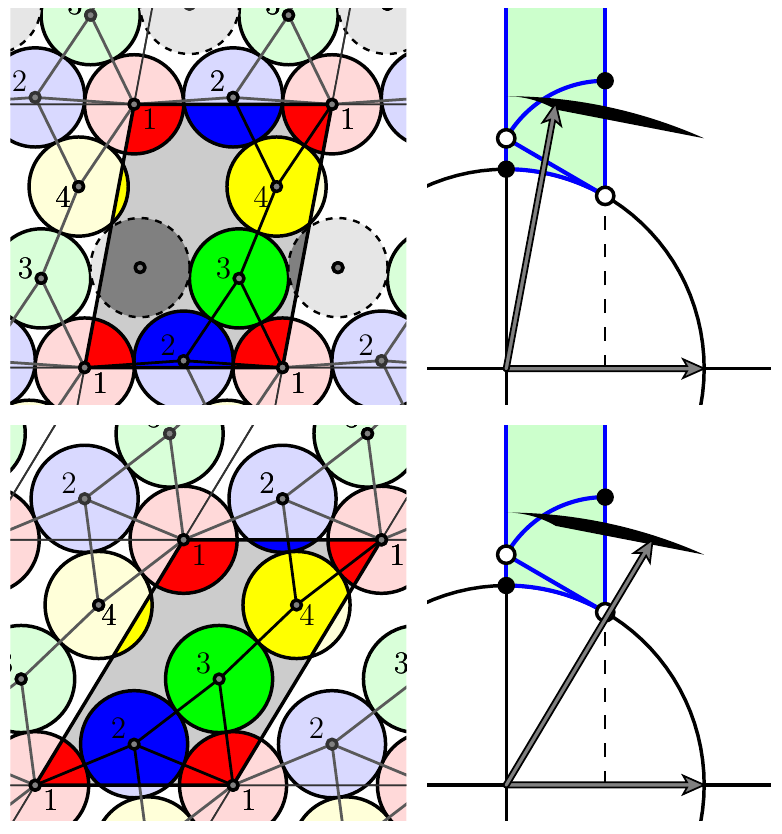}
\caption{A locally and non-globally maximally dense equal circle packing associated to ECG4--2. On the right side, in black, we see the region of the moduli space occupied by these packings. In the top left we can see that it is possible to add a fifth circle (in gray with dashed border) that is not constrained by its neighbors (it is a rattler or free circle). It is not possible to add a sixth circle in the remaining empty space and have the packing be locally maximally dense. In other regions of the moduli space adding a fifth circle is not possible (bottom left).} \label{fig:lmdwfree}
\end{figure}

\begin{figure}
\includegraphics{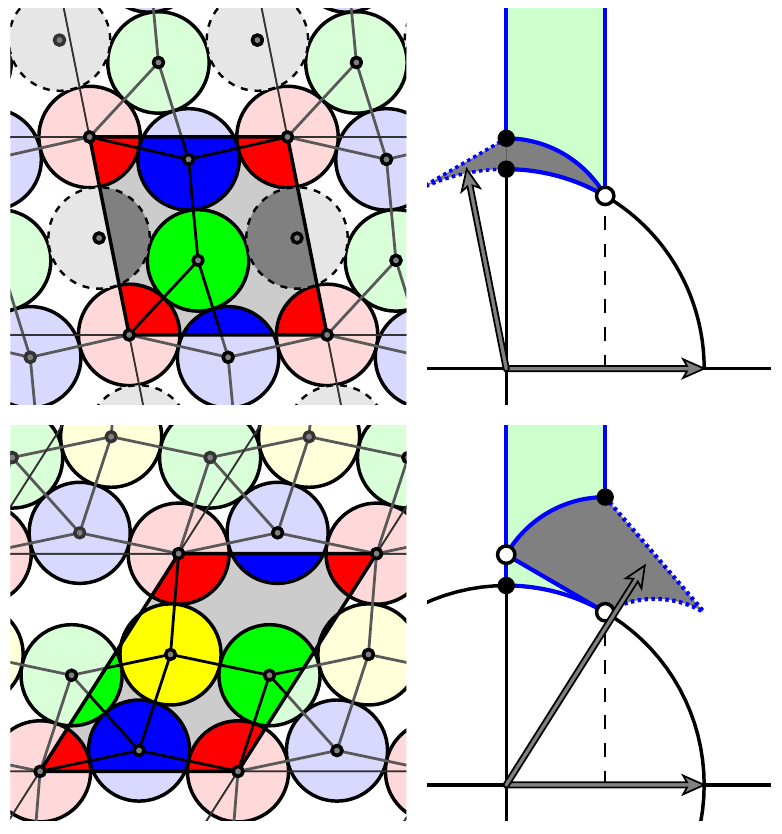}
\caption{The equal circle packings corresponding to ECG1--1 (top row) and ECG9--1 (bottom row)  are globally maximally dense in some regions of the moduli space but in others are only locally maximally dense. Where the gray regions overlap the strip $0\leq x\leq \frac{1}{2}$ they are globally maximally dense, but in the other regions they are only locally maximally dense. For the packing in the upper left, in the non-globally maximally dense region, there is room for another circle (shown in gray with a dashed border). Adding this circle results in a globally maximally dense packing that occupies region $R1_4$ (and four additional tangencies so that the circle is not free).} \label{fig:extended}
\end{figure}

The remainder of the embedded combinatorial graphs embed in a locally maximally dense way on the tori for which they are equal circle packings. First we had to determine which region of the moduli space they occupy. To do this we constrained the embedding so that all lengths were equal and all the angles between adjacent edges were between $\frac{\pi}{3}$ (included) and $\pi$ (excluded).  If an angle is $\pi$ or greater then the packing is not locally maximally dense by Prop.~\ref{prop:VertexLowerBound}. Under these constraints we figured out which region of the moduli space the tori containing the packing occupied.  This is how the regions $Ri_k$ (and boundaries) and the regions in shown in Figures~\ref{fig:lmdwfree} and~\ref{fig:extended} were determined.  To figure out if the equal circle packing was locally maximally dense, we applied  Theorem~6.2 (due to Roth and Whiteley) and Lemma~6.2 from \cite{connelly1}; we found a rigid ordering and a proper stress in the regions of the moduli space indicated. See~\cite{connelly1} for examples and complete details.

The embedded graph for ECG4--2 
corresponds to a locally maximally dense packing that is not globally maximally dense. In Figure~\ref{fig:lmdwfree} you can see the associated equal circle packing and the region of the moduli space of flat tori that it can occupy.  We could transform the region outside of the strip $0 \leq x \leq \frac{1}{2}$ back into this strip in the moduli space, but it would overlap itself and it is more convenient to view it this way.  We note that there is room for a fifth circle that is free in some of the packings. This leads to a locally maximally dense packing of 5 circles that contains a free circle. 

Finally the combinatorial graphs ECG9--1 
and  ECG1--1 
correspond to packings that are locally maximally dense and in some regions of the moduli space are globally maximally dense but not in others. See Figure~\ref{fig:extended}.
 
\section{Existence and Descriptions of the Globally Optimal Packings} \label{sec:Existence}
In this section we discuss the details associated with last part of Step 3 of the overview of the proof found in Subsection~\ref{subsection:overview}.   
We prove the existence of packings that achieve the radii from Theorem~\ref{thm:MainResult} in the regions indicated. To eliminate the translational symmetry we assume that there is always a circle centered at the origin and to simplify the presentation of the coordinates, define the quantity $R_k^{Ri} = \sqrt{16(r_k^{Ri})^2-1}$ where $r_k^{Ri}$ is the expression for the optimal radius for $k$ equal circles in region $Ri_k$. Note that the formulas given in Theorem~\ref{thm:MainResult} for $r_k(x,y)$ are all lower bounded by $\frac{1}{4}$ so $R_k^{Ri}$ is always real and positive.  Let $C_k^{Ri}$ be the list of the centers for $k$ circles in region $Ri_k$ up to equivalence. That is, the coordinates given are in the plane and determine an equivalence class of points in the plane that correspond to the location of the circle centers in the torus.  Note that in the regions $R2_2$, $R3_3$ and $R4_4$ (and all of lower edges of these regions -- See Figure~\ref{fig:ModuliSpaces})  all the circles in the optimal arrangements are self-tangent with radius $\frac{1}{2}$. The arrangements in these regions (except on the lower edge) are far from unique  -- every circle is free. See Section~\ref{sec:selftangent} for complete details.
\subsection{Two Equal Circles} \label{twoequalcircles}
\begin{figure}
\includegraphics{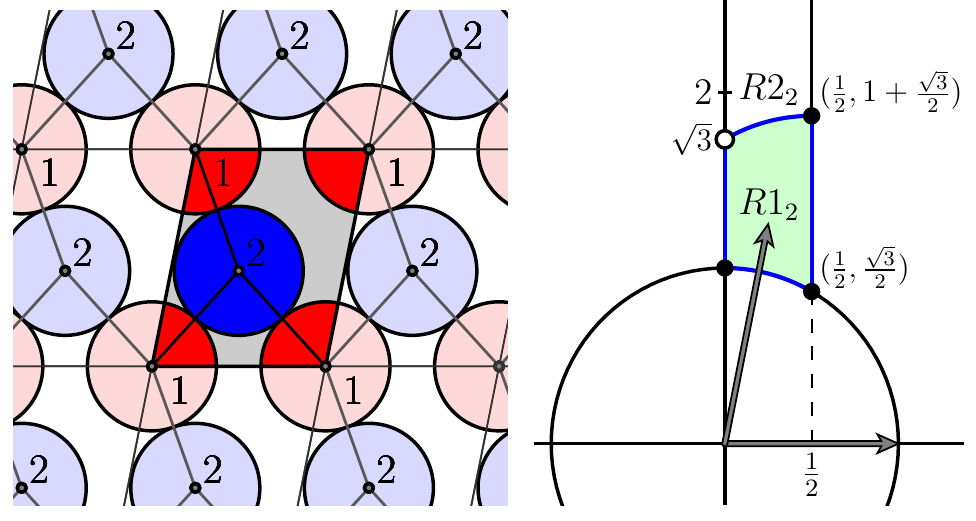}
\caption{On the left is a typical optimal packing of two equal circles in region $R1_2$. On the right is the corresponding location in the moduli space of the torus.} \label{fig:TwoCirclesR1}
\end{figure}

The case of 2 equal circles is included for the sake of completeness. The optimal arrangements and radii are proved in~\cite{przeworski}. However, the tools outlined in this article imply that after fixing the location of the center of the first circle at the origin, the location of the second circle center in the optimal arrangement must be at the circumcenter of the triangle with vertices the origin, $\textbf{v}_1$ and $\textbf{v}_2$.  If this were not the case, then there would be fewer than 3 tangencies and the arrangement could not be locally maximally dense. Note that the bound in Proposition~\ref{prop:GraphProperties} part~(\ref{blah}) is true even for $n=2$.  This implies that
\begin{equation*}
	C_2^{R1} = \bigg\{\ (0,0), \left(\frac{1}{2},\frac{R_2^{R1}}{2}\right) \bigg\}\ .
\end{equation*}
In this region (including the lower edge without the left endpoint and the right edge without the upper endpoint) there are three tangencies. Along the left edge (except for the upper endpoint), when the torus is rectangular, there are four tangencies and along the top edge (excluding the left endpoint) all the circles are self-tangent with radius $\frac{1}{2}$ and there are five tangencies. The only exception to this when the lattice has $\textbf{v}_2 = \langle 0, \sqrt{3} \rangle$ and the packing forms the triangular close packing (each circle is tangent to six others) with six tangencies. See Figure~\ref{fig:TwoCirclesR1} for a typical optimal packing in this region.

\subsection{Three Equal Circles}
For three equal circles, notice that the moduli space is broken into three regions and one can check that the radius is a continuous function in the moduli space. It is interesting to note that on the boundary between regions $R1_3$ and $R2_3$, the radius is actually constant at $\frac{1}{\sqrt{12}}$. To establish the existence of the packing consider the following locations for the circles: 
\begin{align*}
	C_3^{R1} =& \bigg\{\ (0,0), \left(\frac{1}{2},-\frac{R_3^{R1}}{2}\right),\left(\frac{\sqrt{3} R_3^{R1}+1}{4},\frac{\sqrt{3}-R_3^{R1}}{4}\right) \bigg\}\ \\
	C_3^{R2} =& \bigg\{\ (0,0), \left(\frac{1}{2},-\frac{R_3^{R2}}{2}\right),\left(\frac{1}{2}, \frac{R_3^{R2}}{2}\right) \bigg\}\ .
\end{align*}

\begin{figure}
\includegraphics{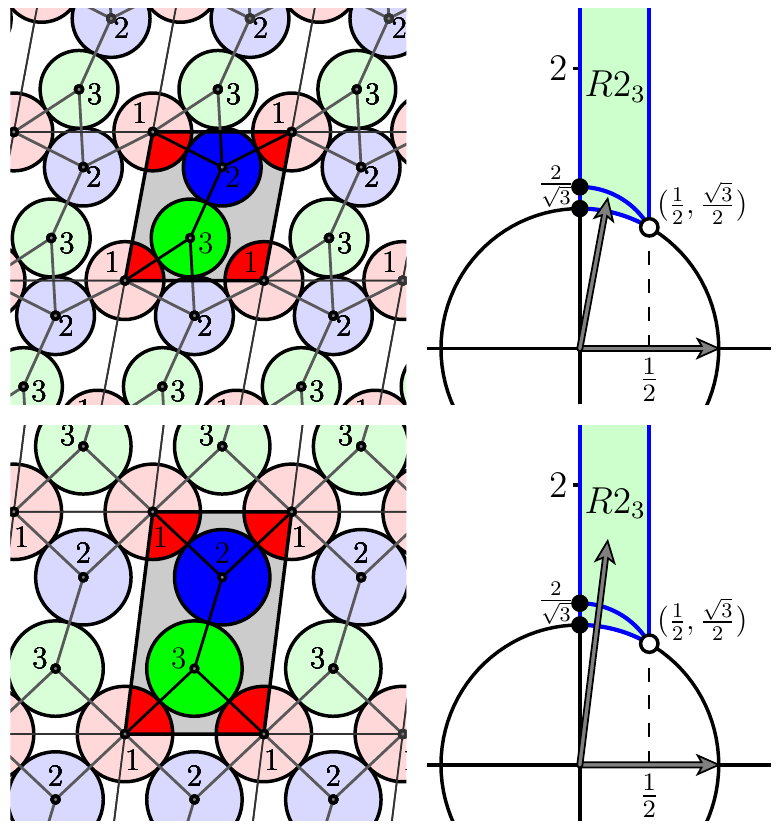} 
\caption{The top row shows a typical optimal packing (associated to ECG1--1) of three equal circles in region $R1_3$ and corresponding location in the moduli space. The bottom row shows the same images only for optimal packings in region $R2_3$ (associated to ECG1--2).} \label{fig:ThreeCirclesR1R2}
\end{figure}

\subsubsection{Region $R1_3$: Three Equal Circles} 
The typical packing in the interior of this region (and including the left and lower edges without the right endpoint) has five tangencies and is associated to ECG1--1. Along the upper edge of this region a sixth tangency is formed so that circle 3 is tangent to circle 1 in two different ways so that, with the edges from circle 2,  two equilateral triangles are formed (associated to ECG 2--1). The only exception to this is when the lattice  is triangular (i.e.  $\textbf{v}_2 = \langle \frac{1}{2},\frac{\sqrt{3}}{2} \rangle$) and the packing forms the triangular close packing with nine tangencies (associated to ECG3--1).  A typical optimal packing from region $R1_3$ is shown the bottom row of Figure~\ref{fig:ThreeCirclesR1R2}.

\subsubsection{Region $R2_3$: Three Equal Circles} 
The typical packing in the interior of this region (and including the left edge without the endpoints) has five tangencies and is associated to ECG1--2.  See the top row of Figure~\ref{fig:ThreeCirclesR1R2}. The right edge (excluding the endpoints) of this region adds a sixth tangency so that the packing graph is a union of rhombi when circle 2 is tangent twice to circle 3 (associated to ECG2--3). Along the upper edge of this region all circles are self tangent with radius $\frac{1}{2}$ and there are eight tangencies except at the righthand endpoint (where $\textbf{v}_2 = \langle \frac{1}{2},\frac{3\sqrt{3}}{2}\rangle$) where the triangular close packing is formed with nine tangencies.  See the left side of Figure~\ref{fig:3} for a typical packing along this circular edge. Note that this packing graph is not on our list in Figures~\ref{fig:toroidalEmbeddings3} or~\ref{fig:toroidalEmbeddings4} because all the circles are self-tangent and the packing graph contains a loop. 

\subsection{Four Equal Circles}
For 4 circles, remarkably, the locations of the centers of circles of the optimal arrangements in regions 1 and 2 is the same relative to the corresponding radius, so the following is true for $i=1$ and $i=2$:
\begin{equation*}
	C_4^{Ri} = \left\{ (0,0), \left(\frac{1}{2},\frac{R_4^{Ri}}{2} \right),\left(\frac{1-\sqrt{3} R_4^{Ri}}{4},\frac{R_4^{Ri}+\sqrt{3}}{4}\right),\left(\frac{3-\sqrt{3} R_4^{Ri}}{4} ,\frac{3 R_4^{Ri}+\sqrt{3}}{4}\right) \right\}. 
\end{equation*}
For region 3, the centers are 
\begin{equation*}
	C_4^{R3} =  \left\{ (0,0), \left(\frac{1}{2}, \frac{R_4^{R3}}{2}\right), \left(0, R_4^{R3}\right),\left(\frac{1}{2}, \frac{3R_4^{R3}}{2}\right) \right\}.
\end{equation*}

\begin{figure}
\includegraphics{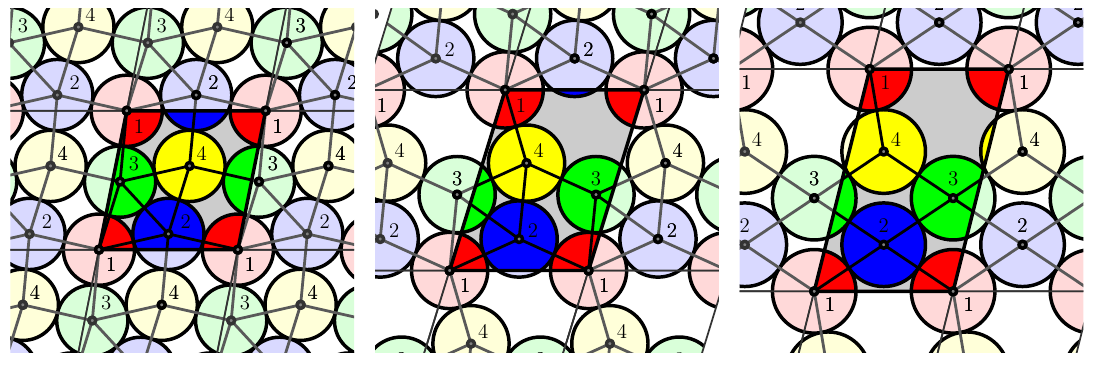}
\caption{On the left (middle/right) is a typical optimal packing of four equal circles in region $R1_4$ ($R2_4$/$R3_4$). These packings are associated to ECG18--1 (ECG9--1/ECG7--1).} \label{fig:FourCirclesR1R2R3}
\end{figure}
\subsubsection{Region $R1_4$: Four Equal Circles} 
The optimal packings in the interior of this region (and including the lower edge without the right endpoint) have nine tangencies and the packing graph is the union of two triangles and three rhombi and is associated to ECG 18--1. Along the left edge (excluding the upper endpoint), circle 4 becomes tangent to the lift of circle 1 at $\textbf{v}_1+\textbf{v}_2$ and the packing graph creates an Archimedean-like tiling of four equilateral triangles and 2 rhombi (associated to ECG20--1). Along the top edge (excluding both endpoints), circle 4 becomes tangent to the lift of circle 1 at $\textbf{v}_2$ and the packing graph creates a different Archimedean-like tiling of four equilateral triangles and 2 rhombi (associated to ECG20--2). The exception to this is when the lattice becomes triangular  (with  $\textbf{v}_2 = \langle \frac{1}{2},\frac{\sqrt{3}}{2} \rangle$ or  $\textbf{v}_2 = \langle 0,\frac{2}{\sqrt{3}} \rangle$) and the packing forms the triangular close packing with twelve tangencies (associated to ECG23--1 or ECG23--2).  A typical optimal packing from region $R1_4$ is shown the left side of Figure~\ref{fig:FourCirclesR1R2R3}. 

\subsubsection{Region $R2_4$: Four Equal Circles}
On the interior of this region (and including the right edge without the lower endpoint) the optimal packings have 8 tangencies and is associated to ECG9--1. Along the top edge (excluding the left endpoint), circle 2 becomes tangent to circle 3 in two different ways and the packing graph is a tiling of 4 triangles and a hexagon (associated to ECG16--1).  It is interesting to note that on the boundary between regions $R2_4$ and $R3_4$, the radius is actually constant at $\frac{1}{\sqrt{12}}$. A typical optimal packing from region $R2_4$ is shown the middle portion of Figure~\ref{fig:FourCirclesR1R2R3}.

\subsubsection{Region $R3_4$: Four Equal Circles} 
The typical optimal packing in the interior of this region (and including the right edge without either endpoint) has seven tangencies and associated to ECG7--1. See the right side of Figure~\ref{fig:FourCirclesR1R2R3}. Along left edge (without either endpoint),  circle 1 becomes tangent to circle 4 in two different ways and the packing graph is the union of four rhombi (associated to ECG 13--1). Along the upper edge of this region all circles are self tangent with radius $\frac{1}{2}$ and there are eleven tangencies except at the lefthand endpoint (where $\textbf{v}_2 = \langle 0,2\sqrt{3} \rangle$) where the triangular close packing is formed with twelve tangencies. See the right side of Figure~\ref{fig:3} for a typical packing along this circular edge. Note that this packing graph is not on our list in Figures~\ref{fig:toroidalEmbeddings3} or~\ref{fig:toroidalEmbeddings4} because all the circles are self-tangent and the packing graph contains a loop.  

In conclusion, we have used rigidity theory to delineate the properties of the combinatorial graphs associated to a locally maximally dense packing with three or four equal circles. There was a list of 23 combinatorial graphs that could possibly be associated to such an optimal equal circle packing. Using Edmond's Permutation Technique we were able to make a list of 103 distinct topologically embedded toroidal graphs.  We were able to eliminate many of these using various techniques and this left us with 15 toroidally embedded graphs that could be associated to an optimal equal circle packing. One of these was locally (and not globally) maximally dense, 2 embedded locally (and not globally) maximally densely on some tori and globally maximally densely on others, and the remaining 12 embedded in a globally maximally dense way. This exhaustive search had to include all globally maximally dense packings (which have to exist because our packing domain is compact) and therefore proves that the radii in Theorem~\ref{thm:MainResult} are globally optimal.

\bibliographystyle{alpha}
\bibliography{ToroidalOptimalEqualCirclePacking}

\end{document}